\documentclass[reqno]{amsart}
\usepackage{fullpage}
\usepackage{amssymb,amscd}

\numberwithin{equation}{section}
\newtheorem{Theorem}[equation]{Theorem}

\newtheorem{Lemma}[equation]{Lemma}
\newtheorem{Proposition}[equation]{Proposition}
\newtheorem{Corollary}[equation]{Corollary}

\theoremstyle{definition}

\newtheorem{Definition}[equation]{Definition}

\theoremstyle{remark}

\newtheorem{Remark}[equation]{Remark}

\newcommand{\eqdef}{\overset{\text{def}}{=}}
\newcommand{\psb}[1]{[ \! [#1] \! ]}

\newcommand{\tensor}[1]{\underset{#1}{\otimes}}

\newcommand{\slot}{\,-\,}

\newcommand{\Z}{\mathbb{Z}}

\newcommand{\R}{\mathbb{R}}
\newcommand{\C}{\mathbb{C}}

\newcommand{\xra}[1]{\xrightarrow{#1}}
\renewcommand{\:}{\colon}

\newcommand{\borel}[1]{{#1}_{\circle}}

\renewcommand{\circle}{\mathbb{T}}
\newcommand{\T}{\circle}

\newcommand{\compact}[1]{{#1}^{+}}
\newcommand{\stalk}[1]{{#1}_{0}}
\newcommand{\stalkata}[1]{{#1}_{a}}
\newcommand{\crep}[1]{\C (#1)}
\newcommand{\cx}[1]{{#1}^{\C}}

\DeclareMathOperator{\Div}{div}

\newcommand{\ET}{E}

\newcommand{\eqvsig}{\Sigma}
\newcommand{\eqvtht}{\Theta}

\newcommand{\h}{\mathfrak{h}}
\newcommand{\HT}{H_{\circle}}

\DeclareMathOperator{\img}{Im}

\newcommand{\I}{\mathcal{I}}
\newcommand{\Izero}[1]{\mathcal{I}_{#1} (0)}

\newcommand{\lchname}{c}
\newcommand{\lch}[1]{c (#1)}
\newcommand{\lm}{\lambda}
\newcommand{\lmh}{\tfrac{\lambda}{2}}

\newcommand{\mero}[1]{\mathcal{M}_{#1}}
\newcommand{\moeight}{MO\langle 8 \rangle}

\newcommand{\nh}{h}
\newcommand{\nhr}{1/h}

\newcommand{\norm}[1]{|#1|}
\newcommand{\Nu}{\mathcal{V}}

\renewcommand{\O}[1]{\mathcal{O}_{#1}}

\newcommand{\phalf}{\tfrac{p_{1}}{2}}

\newcommand{\point}{\ast}
\newcommand{\poth}{PT}

\DeclareMathOperator{\rank}{rank}
\newcommand{\red}[1]{r #1}
\newcommand{\restr}[1]{\rvert_{#1}}

\newcommand{\sg}{\sigma}
\newcommand{\signs}{\delta}

\newcommand{\spin}{\mathrm{Spin}}

\DeclareMathOperator{\Sym}{Sym}

\renewcommand{\th}{\mbox{${}^{\mathrm{th}}$}}
\newcommand{\thalf}{\tfrac{1}{2}}
\newcommand{\trans}{\tau}

\newcommand{\uhalf}{u^{\frac{1}{2}}}
\newcommand{\umhalf}{u^{-\frac{1}{2}}}

\newcommand{\ur}[1]{{#1}^{\R}}

\begin{document}
\title{The Witten genus and equivariant elliptic cohomology}

\author{Matthew Ando}
\email{mando@math.uiuc.edu}
\thanks{The first author was supported NSF grant DMS---0071482}

\author{Maria Basterra}
\email{basterra@math.uiuc.edu}

\address{Department of Mathematics \\
         University of Illinois at Urbana-Champaign \\
         Urbana IL 61801}
 
\date{version 6.4, July 2001}

\begin{abstract}
We construct a Thom class in complex equivariant elliptic cohomology
extending the equivariant Witten genus.  This gives a new proof of the
rigidity of the Witten genus, which exhibits a close relationship to
recent work on non-equivariant orientations of
elliptic spectra. 
\end{abstract}

\maketitle

\section{Introduction} \label{sec-intro}

Let $\T$ denote the circle group.  If $X$ is a $\T$-space, let
$\borel{X}$ denote the Borel construction $E\T\times_{\T} X$.  
Let
$\HT$ denote Borel $\T$-equivariant ordinary cohomology
with complex coefficients.  Choose
a generator of the complex character group of $\T$ and so an isomorphism 
\[
     \HT (pt) \cong \C[z].
\]
Suppose that $\tau$ is a complex number with positive imaginary part.
Consider the lattice $\Lambda=2\pi i \Z + 2 \pi i \tau \Z$, and let
$C=\C/\Lambda$ be the 
associated elliptic curve.  Grojnowski
(\cite{Grojnowski:Ell}; see also \cite{Rosu:Rigidity,GKV:Ell})
has
defined a $\T$-equivariant elliptic cohomology functor $\ET$ from compact
$\T$-spaces to sheaves of algebras over $\O{C}$, equipped with a
canonical isomorphism 
\[
     \stalk{\ET (X)} \cong 
      \HT (X)\otimes_{\C[z]}\stalk{\O{C}}
\]
of stalks at the identity of $C$.  
(For similar constructions 
for $K$-theory see \cite{BBM:ced,DV:ced,BG:echedf,RK:kthy}.) 

The parameter $\tau$ determines a sigma function $\sg = \sg (z,\tau)$;
it is an odd holomorphic function of $z$ which vanishes to first order at
the points of $\Lambda$; see \S\ref{sec-witten-genus}.  Its Taylor
expansion at the origin determines a Thom 
class $\borel{\phi (V)} \in \HT
(V,V_{0})\tensor{\C[z]}\stalk{\O{C}}$  
for any oriented $\T$-vector bundle $V$, which we
call the ``equivariant $\sigma$-orientation''.  The associated genus
is the Witten genus (\cite{AHS:ESWGTC}; see also \S\ref{sec-witten-genus}).

If $V$ is a vector bundle over a compact space $X$, then we define
$\ET (V)$ to be the reduced $E$-cohomology of the Thom space $X^{V}$;
it is a sheaf of $\ET (X)$-modules, called the
Thom sheaf of $V$.    If $V$ is a spin $\T$-vector bundle, then 
(see Proposition \ref{t-pr-thom-sheaf}) $\ET (V)$ is an invertible
sheaf of $\ET (X)$-modules.  More generally, if $W$ is a virtual spin
$\T$-vector bundle over $X$, then there are
genuine spin $\T$-vector bundles $T$ and $V$ such that $W=V-T$, and we may
define $\ET (V-T)$ to be 
the invertible sheaf 
\[
    \ET (W) = \ET (T)^{-1}\otimes_{\ET (X)} \ET (V),
\]
as the right-hand-side is independent of $V$ and $T$ up to canonical
isomorphism (see~\S\ref{sec-thom-sheaf-invertible}).

We recall that if $W$ is a virtual spin vector bundle over $X$, then 
there is a 
characteristic class $\phalf (W) \in H^{4} (X;\Z)$, twice which is the
first Pontrjagin class.   Our main result is the construction of
a Thom class in $\ET$-theory.

\begin{Theorem} \label{t-th-thom-sigma}
Suppose that $\T$ acts non-trivially on a compact smooth manifold
$M$.  Suppose that $W$ is an virtual spin $\T$-vector bundle
over $M$ such that  
\begin{subequations} \label{eq-ccr}
\begin{align}
    w_{2} (\borel{W}) & = 0 \label{eq-swcr}\\
    \phalf (\borel{W}) & = 0. \label{eq-phcr}
\end{align}
\end{subequations}
Then the invertible sheaf of $\ET (X)$-modules $\ET (W)$ 
has a global section $\eqvsig$, whose stalk at the origin is the equivariant
$\sigma$-orientation. 
\end{Theorem}

In the Theorem we have used the fact that the map $V\mapsto
\borel{V}$, which sends a $\T$-vector bundle to its Borel construction,
defines a homomorphism $KO_{T} (X) \rightarrow KO (\borel{X})$.  We
have also used the fact that the Whitney sum formula implies that the
second Stiefel-Whitney class is well-defined for oriented virtual
vector bundles, and the class $\phalf$ is well-defined for virtual
bundles with $w_{2}=0$.  More precisely, we have maps of infinite
loop spaces
\[
\begin{CD}
BSpin @>>> BSO @> w_{2} >> K (\Z/2,2) \\
@V \phalf VV \\
K (\Z,4),
\end{CD}
\]
in which the row is a fibration.

Now consider the case that $W=V-T$, where $T$ is the tangent space of $M$.
If $T$ is spin so that Proposition
\ref{t-pr-thom-sheaf} defines $\ET (-T)$, then the Pontrjagin-Thom
construction provides a map of sheaves 
\[
     \ET (V-T) \xra{\zeta^{*}}\ET(-T) \xra{\poth} \ET (*)\cong \O{C},
\]
where $\zeta$ is the relative zero section.  The class $\stalk{\poth\zeta^{*}
(\eqvsig)}$ is the equivariant Witten genus 
of $M$ twisted by $V$ (see \S\ref{sec-witten-genus}).  Since
$\poth \zeta^{*} (\eqvsig)$ is a global section of $\O{C}$, the twisted Witten
genus is ``rigid'':

\begin{Corollary} \label{t-co-rigidity}
If $T$ is the tangent bundle of $M$, the characteristic classes of $V-T$
satisfy \eqref{eq-ccr},  and $T$  (and so also $V$) is spin, 
then the equivariant Witten genus of $M$ twisted by $V$ is constant.
\end{Corollary}

This result was discovered by Witten, who gave a physical proof in
\cite{Witten:EllQFT}.  The first mathematical proofs were
given by Taubes and Bott-Taubes \cite{Taubes:Rigidity,BottTaubes:Rig}.
Subsequently Kefeng Liu gave a shorter proof and found many new
cases \cite{Liu-0,Liu-1,Liu-2,Liu-3,Liu-4,Liu-5}.

Besides giving  an appealing relationship with equivariant
elliptic cohomology, our construction has the virtue of 
axiomatizing the properties of the $\sg$-function on which it depends.
We construct global sections of $\ET  (W)$ for a class of
theta functions described in 
\S\ref{sec-theta-fns};  Theorem \ref{t-th-thom-sigma} is a special case
of Theorem \ref{t-th-thom-theta}.

We were led to the theta functions of \S\ref{sec-theta-fns} by
the work of Hopkins, Strickland, and the first author.
A theta function such as we consider descends to a ``$\Sigma$
structure'' on the 
elliptic curve $C$, in the sense of \cite{Breen:FonctionsTheta}.
Ando-Hopkins-Strickland have shown that, if $E$ is a (non-equivariant)
elliptic cohomology theory with 
elliptic curve $C$, then maps of ring spectra $\moeight \to E$ are
given by  $\Sigma$ structures on the ideal sheaf $\Izero{\widehat{C}}$ of
the origin in the formal group $\widehat{C}$ of $C$.   

The Weierstrass $\sg$-function given in \S\ref{sec-witten-genus}
descends to the unique $\Sigma$ 
structure on the ideal sheaf $\Izero{C}$, and so gives a
$\Sigma$ structure on $\Izero{\widehat{C}}$ by restriction.  It follows
that the $\sg$-orientation is the
\emph{unique natural} orientation from $\moeight$ to 
elliptic spectra.   Theorem~\ref{t-th-thom-sigma} establishes a
fundamental relationship between these results and the rigidity
theorems for equivariant elliptic genera.  It 
seems reasonable to hope that a (rational) equivariant
analogue of the methods of \cite{AHS:ESWGTC} will produce a functorial
equivariant $\moeight$ 
orientation to the rational equivariant elliptic spectra of Greenlees, 
Hopkins, and Rosu \cite{GHR:Ell}, without the elaborate calculations given 
below.  We will return to that problem in another paper.  

This paper was also inspired by the work of Haynes Miller and Ioanid Rosu
\cite{Miller:Ell,Rosu:Rigidity} .  Miller proposed, before even 
Grojnowski's functor was available, that an uncompleted equivariant
elliptic cohomology should offer a proof of the rigidity theorems.
Rosu and the first author observed that the rigidity
theorems follow from the existence of a global section of the Thom
sheaf in Grojnowski's  theory.  Rosu 
constructed a Thom class whose value at the origin gives the Ochanine
genus.  In fact, he showed that the 
``transfer'' argument of Bott and Taubes is the essential ingredient
in the construction of the Thom section.  
The rigidity of the Ochanine genus and Rosu's analysis require only
that the bundles in question be spin bundles of even rank.
This paper represents a first attempt to understand the Thom classes whose
rigidity require the restrictions \eqref{eq-ccr}.

In fact Theorems \ref{t-th-thom-sigma} and \ref{t-th-thom-theta}
really only require  
\begin{align*}
    w_{2} (\borel{W}) & = 0 \\
    p_{1} (\borel{W}) & = 0,
\end{align*}
as one can see by investigating the use of these conditions in
\S\ref{sec-ccr-consequences}: the vanishing of $\phalf (\borel{W})$ is
used only in rational expressions.  Indeed Liu \cite{Liu-5} has shown that
Corollary \ref{t-co-rigidity} holds without any condition on $w_{2}
(\borel{W})$.   The vanishing of $w_{2} (\borel{W})$ seems to be
necessary for Theorem \ref{t-th-thom-sigma}, and the full conditions
\eqref{eq-ccr}  will certainly be required for any natural equivariant
$\moeight$ orientation to the category of equivariant elliptic
spectra.  

This can already be seen in \S\ref{sec:virtual-reps}, where we
investigate the situation in the case that $W$ is a virtual complex
representation of $\T$ (so $M$ is a point).  In that case too only
the condition $p_{1} (\borel{W}) = 0$ is required to construct a
global section of $\ET (W)$.  However, we do choose a generator $z$ of
the character group of $\T$.  The section of $\ET (W)$ we construct is
independent of this choice only when $w_{2} (\borel{W}) = 0$. 

Our formulation leads to proofs of many of the rigidity theorems in
the literature, including the other genera considered in
\cite{Witten:EllQFT,BottTaubes:Rig} and the twisted loop group
genera of \cite{Liu-0}.  In the interest of brevity, we shall return to
these issues at another time.

\section{Equivariant elliptic cohomology}

In this section, we briefly review Grojnowski's equivariant elliptic
cohomology, following \cite{Grojnowski:Ell} and \cite{Rosu:Rigidity}.
Given a complex elliptic curve $C\cong\C/\Lambda$, we shall describe a 
functor $\ET$ from compact $\T$-spaces to sheaves of
$\Z/2$-graded analytic $\O{C}$-algebras.  The construction uses the
analytic covering 
\[
     \C \rightarrow C
\]
to assemble $\ET (X)$ from the complex $\T$-equivariant cohomology of
$X$.

\subsection{The elliptic curve}

We fix a complex elliptic curve $C$ over $\C$, with an analytic isomorphism 
\[
    C\cong \C/ \Lambda
\]
We write $\pi$ for the covering map 
\[
    \C \xra{\pi} C.
\]
If $V$ is an open set in a  complex analytic variety, then we write
$\O{V}$ for the sheaf of holomorphic functions on $V$, and $\mero{V}$
for the sheaf of meromorphic functions.

Both $\C$ and $C$ are abelian topological groups.  If $G$ is an abelian
topological group and $g\in G$, then we write $\trans_{g}$ for the
translation map; and if $V\subset G$ is an open set, then we write
\[
    V-g \eqdef \trans_{-g} (V).
\]

\begin{Definition} \label{def-small}
An open set  $U$ of 
$C$ is \emph{small} if $\pi^{-1}U$ is a disjoint union of connected
components $V$ such that $\pi\restr{V}: V\to U$ is an isomorphism.
\end{Definition}  

If $U$ is small and $V$ is a component of $\pi^{-1}U$, then 
the covering map induces an isomorphism 
\[
    (\pi\restr{V})_{*}\O{V}\cong \O{U}.
\]
In particular, if $U$ contains the origin of $C$, then there
is a unique component $V$ of $\pi^{-1}U$ containing $0$.  This
determines a $\C[z]$-algebra structure on $\O{U}$, and a
$\C[z,z^{-1}]$-algebra structure on  $\O{U}\restr{U\backslash 0}.$  

\subsection{Adapted open cover of an elliptic curve}

Now suppose that $X$ is a compact $\circle$-manifold.  
If $a$ is a point of $C$, then we define 
\[
    X^{a} = \begin{cases}
    X^{\circle[k]} & a \text{ is of order exactly }k\text{ in }C \\
    X^{\circle} & \text{otherwise}.
\end{cases}
\]
Let $N\geq 1$ be an integer.

\begin{Definition} \label{def-special}
A point $a \in C$ is \emph{special to level }$N$ for $X$ if $X^{Na}\neq
X^{\circle}$.  
\end{Definition}

If $V$ is a $\circle$-bundle over a compact $\circle$-space $X$, then
it is convenient to consider a few additional points to be special.  For each
component $F$ of $X^{\T}$, there are integers $m_{j}$ and an
isomorphism of real $\T$-vector bundles
\[
    V\restr{F} \cong V (0) + \bigoplus \ur{V (m_{j})}.
\]
Here $V (0)$ is the summand of $V\restr{F}$ on which $\T$ acts
trivially, $V (m_{j})$ is a complex vector bundle on which $z\in
\T$ acts by fiberwise by multiplication by the complex number
$z^{m_{j}}$, and if $W$ is a complex vector bundle we write $\ur{W}$
for the underlying real vector bundle.  The $m_{j}$ are called
\emph{exponents} or \emph{rotation numbers} of $V$ at $F$.  Let
$\compact{V}$ denote the one-point compactification of $V$.

\begin{Definition} \label{def-special-vector-bundle}
A point $a$ in $C$ is \emph{special to level} $N$ for $V$ if it is
special for $\compact{V}$ or if for some component $F$ of $X^{T}$ there is a
rotation number $m_{j}$ of $V$ such that $m_{j}Na = 0$.  
\end{Definition}

In any case, for fixed $N$ the set of points which are special to level $N$
is a finite subset of the torsion subgroup of $C$.  Mostly we shall fix an
$N$, and say simply that $a$ is special.

\begin{Definition} \label{def-adapted}
An indexed open cover $\{U_{a} \}_{a \in C}$ of $C$ is \emph{adapted to} $X$
or $V$ (\emph{to level $N$}) if
it satisfies the following.
\begin{enumerate}
\item [1)]  $a$ is contained in $U_{a}$ for all
$a\in C$.  
\item [2)] If $a$ is special and $a \neq b$, then $a\not\in U_{a}\cap U_{b}$.
\item [3)]  If $a$ and $b$ are both special and $a\neq b$, then the
intersection $U_{a}\cap U_{b}$ is empty.
\item [4)]  If $b$ is ordinary, then $U_{a}\cap U_{b}$ is non-empty
for at most one special $a$.
\item [5)]  Each $U_{a}$ is small.
\end{enumerate}
\end{Definition}

\begin{Lemma} \label{t-le-adapted-covers}
$C$ has an adapted open cover, and any two adapted open covers
have a common refinement.  \qed
\end{Lemma}

\subsection{Equivariant cohomology}

We write $\HT$ for Borel $\T$-equivariant ordinary cohomology with
complex coefficients: so if $X$ is a $\T$-space then 
\[
    \HT^{*} (X) = H^{*} (\borel{X}; \C).
\]
We choose a generator of the character group of $\T$, and write $z$
for the resulting generator of $H^{2} (B\T;\Z)$; this gives a
generator of $\HT^{2} (\point)$ which we also call $z$.  We shall
often write $\HT$ for the ring 
$\HT^{*} (\point)\cong\C[z]$.  Moreover we shall consider $\HT$ to be
a subring of the ring $\O{\C} (\C)$ of global analytic functions on
$\C$, by considering $z: \C\to \C$ to be the identity map.

We recall \cite{Quillen:secr} that $\HT$ satisfies a localization
theorem.

\begin{Theorem} \label{t-th-localization}
If $X$ has the homotopy type of a finite $\T$-CW complex (e.g. if $X$
is a compact $\T$-manifold), then the natural map 
\[
     \HT (X) \xra{} \HT (X^{\T})
\]
induces an isomorphism 
\[
 \HT (X)\tensor{\C[z]}\C[z,z^{-1}] \cong 
 \HT (X^{\T})\tensor{\C[z]}\C[z,z^{-1}]. 
\] \qed
\end{Theorem}

\subsection{Elliptic cohomology of a space} 

\label{sec-coh-thom-space}

The equivariant elliptic cohomology of $X$ is a sheaf of
$\O{C}$-algebras over $C$.  We shall describe it as follows.  Let
$\{U_{a} \}_{a\in C}$ be an adapted open cover of $C$.  We first
describe sheaves $\ET (X)_{a}$ over $U_{a}$ for each $a$, and then we
assemble these into a sheaf $\ET (X)$ over all of $C$.

For each $a\in C$, we define an
sheaf of $\O{C}\restr{U_{a}}$-algebras by the formula
\begin{equation}\label{eq-et-local-sheaf}
   \ET (X)_{a} (U) = \HT (X^{a}) \tensor{\C[z]} \O{C} (U-a)
\end{equation}
for $U\subset U_{a}$.  Here $\O{C} (U-a)$ is a $\C[z]$-algebra via
the $\C[z]$-structure on $\O{C}\restr{U_{a}-a}$.  The ring $\O{C} (U)$ acts
by the formula  
\[
    g\cdot (x\otimes y) = x \otimes y \trans_{a}^{*}g.
\]

If $a\neq b$ and $U_{a} \cap U_{b}$ is not empty, then by the definition
\eqref{def-adapted} of an adapted cover, at least one of $U_{a}$ and
$U_{b}$, suppose $U_{b}$, contains no special point.   In particular
we have $X^{b} = X^{\circle}$ and so an isomorphism of $\C[z]$-algebras
\begin{equation}\label{eq-fact-about-X-j}
   \HT (X^{b}) \cong H^{*} (X^{b})\tensor{\C} \C[z].
\end{equation}

\begin{Lemma} \label{t-le-restr-iso}
If $U\subset U_{a}\cap U_{b}$, and $b$ is not special,
then the  inclusion 
\[
    i\: X^{b} \xra{} X^{a}
\]
induces an isomorphism 
\[
\HT (X^{a})\tensor{\C[z]} \O{C} (U-a)  
  \xra{i^{*}\otimes 1}  
\HT (X^{b})\tensor{\C[z]} \O{C} (U-a).
\]
\end{Lemma}

\begin{proof}
If $a$ is not special, then $X^{a}=X^{b}$ and the  result is obvious.
If $a$ is special, then it is not contained in $U$ (by the definition
of an adapted cover), and so $0$ is not contained in $U-a$.    In
particular, $z$ is a unit in $\O{C} (U-a)$.   The localization theorem
\eqref{t-th-localization} gives the result.
\end{proof}

We then define an isomorphism 
\[
    \phi_{ab}\: \ET (X)_{a}\restr{U_{a}\cap
U_{b}} \xra{\cong} 
                \ET (X)_{b}\restr{U_{a}\cap  U_{b}}
\]
of sheaves over $U_{a}\cap U_{b}$ 
as the composition (for $U\subset U_{a}\cap U_{b}$)
\begin{equation}\label{eq-phi-composition}
\begin{split} 
   \HT (X^{a})\tensor{\C[z]} \O{C} (U-a) 
   & \xra{i^{*}\otimes 1}
   \HT (X^{b})\tensor{\C[z]} \O{C} (U-a)  \\
  &\xra{\cong} 
    H^{*} (X^{b})\tensor{\C} \O{C} (U-a)  \\
  & \xra{1\otimes \trans_{b-a} }
    H^{*} (X^{b})\tensor{\C} \O{C} (U-b) \\
  &  \xra{\cong}
    \HT (X^{b})\tensor{\C[z]} \O{C} (U-b).
\end{split}
\end{equation}
The first map is an isomorphism by Lemma \ref{t-le-restr-iso}.  The
rest of the isomorphisms are tautologies, although some use the
isomorphism \eqref{eq-fact-about-X-j}.   

The cocycle condition 
\[
    \phi_{bc} \phi_{ab} = \phi_{ac}
\]
needs to be checked only when two of $a,b,c$ are not special; and in that
case it follows easily from the equation 
\[
     \trans_{c-b} \trans_{b-a} = \trans_{c-a}.
\]
We shall write $\ET ( X )$ for the resulting sheaf over $C$.  One then has
the following.  It was certainly known to Grojnowski
\cite{Grojnowski:Ell}, but for a detailed account the reader may wish
to consult \cite{Rosu:Rigidity}.

\begin{Proposition} \label{t-pr-ell-def}
The sheaf $\ET (X)$ is a sheaf of analytic $\O{C}$-algebras,
which is independent up to canonical isomorphism of the choice of
adapted open cover. \qed
\end{Proposition}

Note that the construction~\eqref{eq-et-local-sheaf} implies that
there is a canonical isomorphism  
\begin{equation} \label{eq-iso-at-stalks}
     \stalkata{\ET (X)} \cong 
      \HT (X^{a})\otimes_{\C[z]}\stalkata{\O{C}}
\end{equation}
of stalks at $a$, as promised in the introduction.

\section{Elliptic cohomology of vector bundles} 
\label{sec-thom-sheaf-invertible}

\subsection{Orientations}
\label{sec:orientations}
Let $HP$ denote \emph{even periodic} cohomology, that is
\[
   HP^{*} (X) =   H^{*} (X;\C[u,u^{-1}]),
\]
where $|u| = -2$.  A power series
\[
      f (z) = z + \text{ higher terms} \in HP^{2} (B\T) \cong \C\psb{z}
\]
satisfying 
\[
     f (-z) = - f (z)
\]
determines an orientation 
\[
    \phi: MSO \rightarrow HP,
\]
characterized by the property that if $L$ is a complex line bundle,
then its Euler class is $f (c_{1}L)$.

If $V$ is  a vector bundle, we write $\phi (V)\in HP (V) = HP (X^{V})$
for the resulting Thom class.  It is \emph{multiplicative}, in the sense that  
\begin{equation}\label{eq-thom-iso-mult}
    \phi (V\hat{\oplus}W) = \phi (V) \wedge \phi (W)
\end{equation}
under the isomorphism 
\[
     (X\times Y)^{V \oplus W} \cong X^{V}\wedge Y^{W}.
\]

\subsection{Equivariant characteristic classes}

If $c$ is a characteristic class of vector bundles and $V$ is a
$\T$-vector bundle over a space $\T$-space $X$, then we write
\[
 \borel{c (V)} \eqdef c (\borel{V}),
\]
provided that makes sense.  We refer to $\borel{c (\slot)}$ as the
``equivariant characteristic class'' associated to $c$.  For example,
if $V$ is oriented, then so is $\borel{V}$, and so
\[
   V\mapsto \borel{w_{2} (V)}
\]
is a characteristic class of oriented $\T$-vector bundles.

If the class $c$ is additive, then we may similarly define $\borel{c
(W)}$ when $W$ is a \emph{virtual} vector bundle.  For example, if $W$
is an oriented virtual $\T$-bundle then we may write
\[
   W = V - V'
\]
where $V$ and $V'$ are oriented $\T$-vector bundles; the quantity
\[
   \borel{w_{2} (W)} = \borel{w_{2} (V)} - \borel{w_{2} (V')};
\]
is easily seen to be independent of the choice of oriented $V$ and
$V'$.  Similarly, if
$W$ is any virtual complex vector bundle, then 
\[
     \borel{c_{1} (W)} = \borel{c_{1} (V)} - \borel{c_{1} (V')}
\]
if 
\[
    W = V - V'
\]
in $K_{\T} (X)$.  

The equivariant Thom isomorphism provides another example.  If
\[
  MSO \xra{\phi} HP 
\]
is an orientation and 
\[
   W = V - V'
\]
where $V$ and $V'$ are oriented $\T$-vector bundles, then we write 
\[
    \borel{\phi (W)} \in HP_{\T} (X^{W}) \cong 
HP_{\T} (X^{V})\tensor{HP_{\T} (X)} HP_{\T} (X^{V'})^{-1}
\]
for the resulting Thom class; and this is independent of the choice of
oriented $V$ and $V'$.

We will use the notation $e_{\phi}$ for the \emph{equivariant} Euler
class associated to $\phi$: so if $V$ is an oriented $\T$-vector
bundle over $X$ then 
\[
   e_{\phi} (V) = \zeta^{*} (\borel{\phi (V)}) \in HP_{\T} (X),
\]
where 
\[
   \zeta: X \rightarrow X^{V}
\]
is the zero section.

\subsection{Analytic orientations}

Let 
\[
\phi: MSO \rightarrow HP
\]
be an orientation determined by a power series $f$ as in
\S~\ref{sec:orientations}.  

\begin{Definition}\label{def-analytic-or} 
The orientation $\phi$ is \emph{analytic}
if $f$ is contained in the subring
$\stalk{\O{\C}} \subset \C\psb{z}$ of germs of holomorphic functions
at $0$; equivalently, if there is a neighborhood $U$ of $0$ in $\C$ on
which the power series $f$ converges to a holomorphic function.  
\end{Definition}
The basic fact about analytic orientations is the following. 

\begin{Proposition}[\cite{Rosu:Rigidity}]
\label{t-pr-analytic-euler-classes} 
If $\phi$ is analytic and $V$ is an oriented $\T$-vector bundle over a
compact $\T$-space $X$, then the equivariant Euler class $e_{\phi}$
associated to $\phi$ satisfies 
\[
   e_{\phi}(V) \in \HT (X)\tensor{\C[z]}\stalk{\O{\C}}.
\]\qed
\end{Proposition}

If $\Phi$ denotes the standard Thom isomorphism, then 
\[
   \borel{\phi (V)} = \frac{e_{\phi} (V)}{e_{\Phi} (V)} \borel{\Phi (V)},
\]
and the ratio of Euler classes is a unit in $\HT
(X)\tensor{\C[z]}\stalk{\O{\C}}.$ 
Of course multiplication by $\borel{\Phi (V)}$ induces an isomorphism 
\[
    \HT (X) \cong \HT (V),
\]
and so we have the following

\begin{Corollary} \label{t-cor-analytic-thom}
There is a neighborhood $U$ of the origin in $\C$ such
that 
\[
    \borel{\phi (V)} \in \HT (V)\tensor{\C[z]} \O{\C} (U), 
\]
and such that multiplication by this class induces an isomorphism of
sheaves 
\[
   \HT (X)\tensor{\C[z]}\O{U} \xrightarrow[\cong]{\phi}    
   \HT (V)\tensor{\C[z]} \O{U}.
\]
In other words, for every open set $U'\subseteq U$, multiplication by
$\borel{\phi (V)}$ induces an isomorphism 
\[
   \HT (X)\tensor{\C[z]}\O{\C} (U') \xrightarrow[\cong]{\phi}    
   \HT (V)\tensor{\C[z]} \O{\C} (U').
\] \qed
\end{Corollary}

\begin{Definition} \label{def-adapted-thom}
Let $V$ be a $\circle$-vector bundle over a compact $\circle$-manifold
$X$, and let $\phi$ be a multiplicative analytic orientation.  An
indexed open cover $\{U_{a} \}_{a\in C}$ of $C$ is 
\emph{adapted to the pair} $(V,\phi)$ to level $N$ if it is
adapted to $V$ to level $N$ (see Definition~\ref{def-adapted}), and 
if for every point $a\in C$, the equivariant Thom class $\phi (V^{a})$
induces an isomorphism 
\[
\HT (X)\tensor{\C[z]}\O{U_{a}-a} \xrightarrow[\cong]{\phi}
\HT (V^{a})\tensor{\C[z]} \O{U_{a}-a}.
\]
\end{Definition}

\subsection{Cohomology of the Thom space}

Now suppose that $V$ is a $\circle$ vector bundle over $X$.  Let
$\ET (V)$ denote the reduced equivariant elliptic cohomology of the
one-point compactification of $V$.  In the case that $V$ is spin,
we give an explicit cocycle which exhibits $\ET (V)$ as an
invertible sheaf of $\ET (X)$-modules.  The main tool which makes this
possible is the following.

\begin{Lemma} \label{t-le-oriented}
Let $V$ be a $\spin$ $\T$-vector bundle over $X$.  For
all $n$, the fixed bundle $V^{\T[n]}$ over $X^{\T[n]}$ is orientable
(that is, it has a Thom isomorphism in ordinary cohomology with integer
coefficients). 
\end{Lemma}

\begin{proof}
If the rank $V$ is even, this is Lemma 10.1 in \cite{BottTaubes:Rig}.
If the rank of $V$ is odd, then we may apply that lemma to $V'=V\oplus
\R$.
\end{proof}

Suppose then that $V$ is a Spin $\T$-vector bundle over $X$.  Let $\phi$ be a multiplicative analytic orientation.  Let $e$ be the
associated  (equivariant) Euler class.  We define a class $[\phi,V]
\in H^{1} (C; \ET(X)^{\times})$ as follows.  
Choose an open cover
adapted to the pair $(V,\phi)$.  
Suppose that $a, b$ are two points of 
$C$, such that $U_{a}\cap U_{b}$  is non-empty: we may suppose that
$b$ is not special.  Suppose that $U\subseteq U_{a}\cap U_{b}$.
Recall that there is an isomorphism 
\eqref{eq-phi-composition} 
\[
   \ET (X) (U)\cong \HT (X^{b})\tensor{\C[z]} \O{C} (U-a).
\]

For each special point
$a$ of $C$, choose an orientation on $V^{a}$ and on $V^{\T}$; we may
do so by Lemma \ref{t-le-oriented}.  This gives equivariant euler classes $e
(V^{a})$ for $a\in C$.

\begin{Lemma} \label{t-le-cocycle-defined}  There is a unit $e (a,b)\in \ET (X) (U)$ such that $e (V^{b})  e (a,b)
= e(V^{a}\restr{X^{b}})$.  Moreover we have 
\[
  e (a,a) = 1
\]
for all $a$, and 
\begin{equation}\label{eq-cocycle-h-1}
  e (a,b)e (b,c) = e (a,c) 
\end{equation}
whenever that makes sense.
\end{Lemma}

\begin{proof}
Without loss of generality we may suppose that $V^{b} = V^{\circle}$.
If in addition $V^{a}= V^{\circle}$, then $e (a,b) = 1$. 
Otherwise, we have $a\not\in U$ so $0\not\in U-a$, and so $z$ is a unit in
$\O{C} (U-a)$.  On each component $F$ of $X^{b}$, there are integers
$m_{j}\neq 0$ and complex vector bundles $V(m_{j})$ over $F$ such that
\[
     V^{a}\restr{F} = V^{b}\restr{F} \oplus \bigoplus_{m_{j}} \ur{V (m_{j})}.
\]
Here $\T$ acts on $V (m_{j})$ fiberwise by the character $u\mapsto
u^{m_{j}}$.  

Let $f (z) = z + \text{higher terms}$ be the characteristic series of
the orientation $\phi$, so that 
\[
     e (L) = f (c_{1} L)
\]
for $L$ a complex line bundle.  Let $x_{j,1},\dots ,x_{j,d_{j}}$ be
the roots of the total Chern class of $V (m_{j})$.  Since $F$ is
compact, the $x_{i,j}$ are nilpotent.  It suffices to check that $e
(a,b)\restr{F}$ is a unit modulo nilpotents in $H (F)\tensor{\C[z]}
\O{C} (U-a)$.  We have 
\begin{align*}
   e (a,b)\restr{F} & = \prod_{j}\prod_{i} f (m_{j}z + x_{j,i})\\
           & \equiv  \prod_{j}\prod_{i} f (m_{j}x) \\
           & \equiv  (\prod_{j} m_{j}^{d_{j}}) z^{\sum d_{j}} (1 +
\text{higher terms in $z$}),
\end{align*}  
where the equivalences are modulo nilpotents.  The result
follows, since the $m_{j}$ are non-zero and $z$ is a unit in $\O{C}
(U-a)$.  

The cocycle condition \eqref{eq-cocycle-h-1} is easy, because as usual
the equation needs only to be verified when at most one of $a$, $b$,
and $c$ is special.
\end{proof}

Let $[\phi,V]\in H^{1} (C;\ET (X)^{\times})$ be the cohomology class
defined by the $e (a,b)$.  Let $\ET (X)^{[\phi,V]}$ denote the
resulting invertible sheaf of $\ET (X)$-modules over $C.$
Explicitly, the sheaf $\ET (X)^{[\phi,V]}$ is assembled from
the sheaves 
\[
    \ET (X)_{a} (U) = \HT (X^{a})\tensor{\C[z]} \O{C} (U-a)
\]
over $U_{a}$, using the sequence of isomorphisms

\begin{equation} \label{eq-twisted-sheaf}
\begin{split}
    \ET_{a} (X) (U) & \cong 
    \HT (X^{a})\tensor{\C[z]} \O{C} (U-a)  \\
      & \xra{i^{*}\otimes 1}
    \HT (X^{b})\tensor{\C[z]} \O{C} (U-a) \\
      & \xra{e (a,b)\cdot}
    \HT (X^{b}) \tensor{\C[z]} \O{C} (U-a) \\
      & \cong 
    H^{*} (X^{b})\tensor{\C} \O{C} (U-a)  \\
      & \xra{1\otimes \trans_{b-a}^{*}} 
    H^{*} (X^{b}) \tensor{\C} \O{C} (U-b) \\
      & \cong
    \HT (X^{b}) \tensor{\C[z]} \O{C} (U-b) \\
      & \cong \ET_{b} (X) (U)
\end{split}
\end{equation}
as transition functions to assemble a sheaf over $C$.

It will be important in \S\ref{sec-thom-class} that we 
allowed the orientation of $V^{\T}$ to vary with the 
special point $a$.  In fact, the resulting sheaf $\ET (X)^{[\phi,V]}$
is independent of the choices, up to canonical isomorphism.

\begin{Proposition} \label{t-pr-thom-sheaf}
If $V$ is a $spin$ $\T$-bundle, then the Thom isomorphism
$\phi$ induces an isomorphism 
\[
     \ET (X)^{[\phi,V]} \cong \ET (V)
\]
of sheaves of $\ET (X)$-modules.
\end{Proposition}

\begin{proof}
Choose a cover adapted to $X$ and the pair $(V,\phi)$.  Suppose that $U\subset
U_{a}\cap U_{b}$, $a$ is special, and $b$ is not.  Then the diagram
\[
\begin{CD}
   \HT (X^{a})\tensor{\C[z]} \O{C} (U-a) 
   @> \phi >> 
   \HT (V^{a}) \tensor{\C[z]} \O{C} (U-a) \\
   @V e (a,b)i^{*} VV @VV i^{*} V \\
   \HT (X^{b})\tensor{\C[z]}\O{C} (U-a)
   @> \phi >> 
   \HT (V^{b}) \tensor{\C[z]} \O{C} (U-a) \\ 
   @V 1 \otimes \trans_{b-a}^{*} VV @VV 1 \otimes \trans_{b-a}^{*} V \\
   \HT (X^{b})\tensor{\C[z]}\O{C} (U-b) @> \phi >> 
   \HT (V^{b}) \tensor{\C[z]} \O{C} (U-b)
\end{CD}
\]
commutes (all the arrows are isomorphisms).  The left column describes
the sheaf $\ET (X)^{[\phi,V]}$, while the right column describes $\ET
(V)$. 
\end{proof}

Now suppose that $W$ is a \emph{virtual} $\T$-bundle.  We may
write 
\begin{equation}\label{eq-W-diff}
   W = V - T
\end{equation}
with $V$ and $T$ genuine $\T$-bundles of even rank.  We also clearly
have 
\[
    W = (V+3T) - 4T, 
\]
and it is easy to check that
\begin{align*}
    w_{1} (4T)_{\T} & = 0 \\
    w_{2} (4T)_{\T} & = 0 \\
    w_{1} (V+3T)_{\T} & = w_{1} (W)_{\T} \\
    w_{2} (V+3T)_{\T} & = w_{2} (W)_{\T}.
\end{align*}

Thus if $W$  is a virtual spin bundle, then we may require
that $V$ and $T$ are spin bundles of even rank.  If
\[
   W = V - T = V' - T'
\]
with $V$, $V'$, $T$, and $T'$ spin bundles of even rank, then 
\[
   V = V' \pm D
\]
where $D$ is a spin bundle.  It follows that we may 
define $\ET (W)$ to be the invertible sheaf
\[
   \ET (W) = \ET (T)^{-1}\tensor{\ET ( X)} \ET (V).
\]
As above, we can as construct a class $[\phi,V-T]\in H^{1} (C,\ET
(X)^{\times})$, equipped with a canonical isomorphism 
\[
    \ET (W)\cong \ET (X)^{[\phi,V-T]}.
\]


\section{Theta functions}
\label{sec-theta-fns}

The orientations in this paper will arise as the Taylor expansions of
certain meromorphic $\theta$ functions.  The $\theta$ functions we
consider are a special case of those considered
in~\cite[Ch. I]{Mumford:AbelianVarieties}, except that Mumford
considers only holomorphic functions.  

Let $\Lambda\subset \C$ be a lattice.
By a \emph{theta function} for $\Lambda$ we mean a meromorphic function
$\theta\: 
\C\to \C$ with the following properties.  First, there is an integer
$N$ such that the zeroes and poles of 
$\theta$ are contained in $N^{-1}\Lambda$.  Second, there are
functions 
\begin{align*}
    \gamma:  \Lambda & \to \C \\
    \lchname:  \Lambda & \to \{\pm 1 \}
\end{align*}
such that, for all $\lambda\in \Lambda$ and $z\in \C$, we have 
\begin{subequations}
\label{eq-gp-theta}
\begin{align}\label{eq-theta-odd} 
\theta (-z) & = - \theta (z)\\
\theta (0)  & = 0 \label{eq-theta-vanishes}\\
\theta' (0) & \neq 0 \label{eq-theta-first-order} \\
 \theta (z + \lambda) & = \lch{\lm} e^{\gamma (\lambda) (z+\lmh)} \theta (z).
 \label{eq-theta-translates-by-char}
\end{align}
\end{subequations}

Using \eqref{eq-theta-translates-by-char} in various ways to
compare $\theta (z+\lambda + \lambda')$ with $\theta (z)$, one finds 
that $\gamma$ is necessarily linear, and that
\[
   \frac{\lch{\lambda + \lambda'}}
        {\lch{\lambda}\lch{\lambda'}} = 
e^{\tfrac{1}{2} (\gamma (\lambda)\lambda' - \lambda \gamma
(\lambda') )}
\]
for $\lambda, \lambda'\in \Lambda$.  It follows that  $\gamma$ must satisfy the
``period relation''
\[
    \gamma (\lambda)\lambda' - \lambda \gamma (\lambda') \in 2 \pi i
\Z.
\]
The linearity of $\gamma$ gives the useful formula
\begin{equation}\label{eq-trans-form-sg-it}
   \theta (z+ \ell \lm) = \lch{\ell \lm} e^{\gamma (\lambda) ( \ell z +
\ell^{2}\lmh)} \theta (z).
\end{equation}

\begin{Remark} \label{rem-theta}
Equation \eqref{eq-theta-translates-by-char} means that $\theta$ is a
meromorphic section of the line bundle 
\[
   \mathcal{L} (\gamma,c) = 
   \frac{\C\times \C}
        {(z,v) \sim (z + \lm, v
      c (\lm) e^{\gamma (\lm) (z + \tfrac{\lm}{2})}), \lm \in \Lambda}
\]
over $\C/\Lambda$, so the theta functions considered here are a
special case of \cite[Ch. I]{Mumford:AbelianVarieties}.
\end{Remark}

The equations (\ref{eq-theta-odd}---\ref{eq-theta-first-order})
imply that the Taylor expansion $\hat{\theta}$ of $\theta$ at the
origin determines a 
multiplicative analytic orientation (see Definition~\ref{def-analytic-or})
\begin{equation}\label{eq-orientation-from-theta}
    \phi: MSO\to HP.
\end{equation}

The genus associated to this orientation is 
\begin{equation}\label{eq-genus-integration-formula}
    M \mapsto \int_{M} \prod_{j=1}^{d} \frac{x_{j}}{\hat{\theta} (x_{j})},
\end{equation}
where $M$ is an oriented  manifold of dimension $2d$, and the $x_{j}^{2}$
are the roots of the total Pontrjagin class of $M$.

\section{Examples}
\subsection{The Witten genus}
\label{sec-witten-genus}

One example of a theta function satisfying \eqref{eq-gp-theta} 
is the Weierstrass $\sg$-function.  We describe  a variant  associated
to the Witten genus.

Let $\sg$ denote the expression
\begin{equation}\label{eq-sg-prod-expr}
    \sg = (\uhalf - \umhalf) 
              \prod_{n\geq 1} \frac{(1-q^{n}u) (1-q^{n}u^{-1})}
                                   {(1-q^{n})^{2}}.
\end{equation}
This may be considered as an element of
$\Z\psb{q}[u^{\pm\frac{1}{2}}]$ which is a  holomorphic function of
$(\uhalf,q)\in \C^{\times}\times D$, where $D = \{q\in \C |
0<\norm{q}<1 \}$. Let $\h = \{\tau\in \C | \img \tau > 0 \}$ be the
open upper half 
plane.  We may consider $\sg$ as a holomorphic function of
$(z,\tau)\in \C\times \h$ by setting 
\begin{align*}
      \uhalf & = e^{\frac{z}{2}} \\
      q  & = e^{2\pi i \tau}.
\end{align*}

For our purposes, it is 
sufficient to fix $\tau$, and consider $\sg$ as a function of $z$
alone.  From the formula \eqref{eq-sg-prod-expr} it 
follows that 
$\sg (z)$ is a theta function of the form \eqref{eq-gp-theta} for the
lattice $\Lambda = 2 \pi i \Z + 2 \pi i \tau \Z$, with 
\begin{align*}
    \gamma (2\pi i a + 2 \pi i b \tau ) & = - b \\
    \lch{\lm} & = \begin{cases}
1 & \lm \in 2\Lambda \\
-1& \text{otherwise.}
\end{cases}
\end{align*}

If $V$ is a complex vector bundle, let
$\red{V} = \rank V - V$ denote the associated virtual bundle of rank $0$.  Let
$S^{k} V$ denote the  $k$ symmetric  
power of $V$, and let 
\[
   \Sym_{t} (V) = \sum_{k\geq 0} t^{k} S^{k}V.
\]
This extends to an operation 
\[
    K (X)\to (1 + t K (X)\psb{t})^{\times}
\]
because of the formula 
\[
\Sym_{t} (V\oplus W) = \Sym_{t} V \cdot \Sym_{t} W.
\]
The genus associated to the (Taylor series expansion of) 
\[
   a (x) = e^{\tfrac{x}{2}} - e^{-\tfrac{x}{2}}
\]
is called the $\widehat{A}$-genus.  
Formulae \eqref{eq-genus-integration-formula} and
\eqref{eq-sg-prod-expr} show that the genus  
associated to $\sg$ is 
\[
     M\mapsto 
\widehat{A} (M; \bigotimes_{n\geq 1} \Sym_{q^{n}} (\red{\cx{TM}}))
\]
(this is explained in~\cite{AHS:ESWGTC}), which is equivalent to formula (27) in \cite{Witten:EllQFT}.

\begin{Remark} 
Let $C$ be the complex elliptic curve $\C/\Lambda$.
The equations \eqref{eq-gp-theta} in this case are descent data for the
``$\Sigma$-structure'' on the ideal sheaf $\Izero{C}$ of the origin in
$C$, in the sense of Breen \cite{Breen:FonctionsTheta}; see
particularly {\S}3.12. 
\end{Remark}

\begin{Remark} \label{rem-relation-classical-sigma}
The sigma function \eqref{eq-sg-prod-expr} is related to the
classical Weierstrass sigma function  
$\sg_{\mathrm{Weierstrass}}$ in for example
\cite{Silverman:EllipticCurves} by the formula 
\[
    \sg_{\mathrm{Weierstrass}} (z) = e^{a z^{2}} \sg (z),
\]
where $a$ is a constant.  It is not hard to check that these define
the same invariants of $\moeight$-manifolds.  The $\sg$-function
\eqref{eq-sg-prod-expr} also 
arises as the $p$-adic $\sg$-function of the Tate curve
\cite{MazurTate:sigma}. 
\end{Remark}

\subsection{The Ochanine genus}\label{sec-ochanine}

Let $C\cong \C/\Lambda$ be an elliptic curve, and let $p$ be a
point of exact order $2$ of $C$.  Let $r_{1}$ and $r_{2}$ be the other
two non-zero points of order two, and let $P,R_{1},R_{2}$ be
representatives of these points, such that $R_{1}+R_{2}=P$.  There is
a meromorphic function $f$  on $C$ with divisor  
\[
       \Div f = (0) + (p) - (r_{1}) - (r_{2});
\]
$f$ is uniquely determined by requiring that its pull-back $s$ to $\C$
have a Taylor series expansion of the form 
\[
   s (z) = z + o (z^{2}).
\]
In fact in this case $s$ is given by the formula 
\begin{equation}\label{eq-s-theta}
    s (z) = \frac{\sg (z)\sg (-R_{1})\sg (-R_{2})\sg (z-P)}
                 {\sg (z-R_{1})\sg (z-R_{2})\sg (-P)}.
\end{equation}
It is odd, and the resulting genus is the elliptic
genus of Ochanine \cite{Ochanine:EllipticGenera} for the pair
consisting of the lattice $\Lambda$ and the point $p$.  
It is customary to consider $s$ as a theta
function for the lattice 
\[
\Lambda' = \{\lm + n P | \lm \in \Lambda, n \in \Z\}.
\]
The function $s$ is not quite periodic with respect to this lattice,
but satisfies 
\[
  s (z + \lm + n P ) = (-1)^{n} s (z)
\]
for $\lm\in \Lambda$.  Thus 
\begin{align*}
     \gamma (\lm + n P ) & = 0 \\
     \lch{\lm + n P } & = (-1)^{n}.
\end{align*}
Explicitly, if $\Lambda=2\pi i \Z + 4 \pi i \tau\Z$ and $P=2\pi i
\tau$, then using \eqref{eq-sg-prod-expr} in \eqref{eq-s-theta} gives
the formula
\[
    s (z) = -2 \frac{1-u}{1+u}
            \prod_{n\geq 1} 
             \frac{(1-q^{n}u) (1-q^{n}u^{-1}) (1+q^{n})^{2}}
                  {(1 + q^{n}u) (1+q^{n}u^{-1}) (1-q^{n})^{2}}.
\]
Theorem \ref{t-th-thom-theta} for this function is due to Rosu
\cite{Rosu:Rigidity}.

\section{The equivariant Thom class}

\label{sec-thom-class}

We fix a meromorphic function $\theta$ satisfying the conditions of
\S\ref{sec-theta-fns}.  We recall~\eqref{eq-orientation-from-theta}
that the Taylor expansion $\hat{\theta}$ of $\theta$ at $0$ determines
a multiplicative analytic orientation 
\[
    MSO \rightarrow HP,
\]
and we write $\phi$ for the resulting Thom isomorphism.

Suppose that $X$ is a compact smooth $\T$-manifold, and that $W$
is a virtual oriented $\T$-bundle over $X$.  The Thom isomorphism $\phi$
gives a generator
$\borel{\phi (W)}$ of $\HT (W)$, 
and so a generator of
\[
   \HT (W) \tensor{\C[z]}\stalk{\O{C}}
\]
which we shall also denote $\borel{\phi (W)}$.

\begin{Theorem}  \label{t-th-thom-theta}
If $W$ is an oriented $\T$-vector bundle over
$X$, and either
\begin{enumerate}
\item $w_{2} (W) = 0$ and the function $\gamma$ of
\eqref{eq-theta-translates-by-char} is identically zero; or
\item the  equations \eqref{eq-ccr} hold, 
\end{enumerate}
then the invertible sheaf of $\ET (X)$-modules $\ET (W)$ has a 
global section $\eqvtht$, such that
\[
    \eqvtht_{0} = \borel{\phi(W)}
\]
under the isomorphism~\eqref{eq-iso-at-stalks} 
\[
   \HT (W) \tensor{\C[z]}\stalk{\O{C}} \cong \stalk{\ET (W)}.
\]
\end{Theorem}

The proof will occupy the rest of this section.  To give it we
fix a cover of $C$ adapted to level $N$ to $(V,\phi)$ and $(T,\phi)$,
where $N$ is large enough that the zeroes and poles of $\theta$ are
contained in $N^{-1}\Lambda$.  We also write 
\[
   W = V - T,
\]
where $V$ and $T$ are spin $\T$-vector bundles of even rank, as in
\S\ref{sec-thom-sheaf-invertible}.  

Let us indicate precisely what it is we must construct.   Since $\phi$
is multiplicative, we have 
\[
\borel{\phi (W)} = \borel{\phi (T)}^{-1}\otimes \borel{\phi (V)}.
\]
Proposition \ref{t-pr-thom-sheaf} shows that it is
equivalent to construct a global section of the sheaf $\ET
(X)^{[\phi,V-T]}$, whose value in 
\[
\ET (X)^{[\phi,V-T]}_{0}\cong \HT (X)\tensor{\C[z]}\O{C} (U_{0})
\]
is $1$.  The formula \eqref{eq-twisted-sheaf} for this
sheaf shows that such a  global section is assembled from sections
$\eqvtht_{a}\in \ET (X)_{a} (U_{a})$ which satisfy the formula 
\begin{equation}\label{eq-mu-transition-fn}
       \eqvtht_{b}   =      \trans_{b-a}^{*}  ( e (a,b) i^{*}\eqvtht_{a}).
\end{equation}
on $U\subset U_{a}\cap U_{b}$.  Because we are using an adapted open
cover, it will suffice in \eqref{eq-mu-transition-fn}  to suppose
that $b$ is not special.  

In order to give the formula for $\eqvtht$, we introduce some notation
and results in section \S~\ref{sec:notations-lemmas}.  The
construction of $\Theta$ using these 
results begins in \S~\ref{sec:constr-thet-ordin}. 

\subsection{Notations and lemmas}
\label{sec:notations-lemmas}

Suppose that $a$ is a special point of order $n$.  
The structure of the formula for $\eqvtht_{a}$ depends
on the parity of $n$, but the two cases share many components.  Let
$\nh=n/2$.  In the following, terms which involve $\nh$ are simply
absent in the case that $n$ is odd.

Let $P$ be a component of the submanifold $X^{\T[n]}$ of points fixed
by $\T[n]$.  Let $F\subset P$ be a component of $X^{\T}$.  We
have decompositions of real $\T$-vector bundles
\begin{align}\label{eq-T-T-sub-zero-to-h}
      T\restr{P} & = T_{0} \oplus T_{\nh} 
                      \oplus \bigoplus_{0<r<\nh} \ur{T_{r}} \\
      V\restr{P} & = V_{0}  \oplus V_{\nh} 
                      \oplus \bigoplus_{0<r<\nh} \ur{V_{r}}. \notag
\end{align}
For $0<r<\nh$, $T_{r}$
and $V_{r}$ are complex vector bundles over $P$ on which $\zeta\in \T[n]$
acts by multiplication by $\zeta^{r}$.  The bundles $T_{0}$ and $V_{0}$
are the summands of $T$ and $V$ on which $\T[n]$ acts trivially; and
if $n$ is even then
$T_{\nh}$ and $V_{\nh}$ are the summands on which $\T[n]$ acts by the
sign representation. 

Over the submanifold $F$ there are decompositions 
\begin{align}
      T_{0}\restr{F} & = 
               T (0) \oplus 
               \bigoplus_{0\neq m_{j}\equiv 0}  \ur{T (m_{j})} 
   \label{eq-T-zero-T-sub-zero}\\
      T_{h}\restr{F} & = 
               \bigoplus_{m_{j}\equiv \nh}  \ur{T (m_{j})} 
  && n \text{ even } 
    \label{eq-T-h-decomp} \\
      T_{r}\restr{F} & = \bigoplus_{m_{j} \equiv r}
                                     \ur{T (m_{j})}   && 0<r<\nh \notag     \\
      V_{0}\restr{F} & = 
                V (0) \oplus 
                \bigoplus_{0\neq m_{j}'\equiv 0}  \ur{V (m_{j}')} 
                                        \notag \\
      V_{h}\restr{F} & = 
               \bigoplus_{m_{j}\equiv \nh}  \ur{V (m_{j})} 
  && n \text{ even }                    \notag  \\
      V_{r}\restr{F} & = \bigoplus_{m_{j}'\equiv r}
                                     \ur{T (m_{j}')}   && 0<r< \nh; \notag
\end{align}
the equivalences are modulo $n$.   For $m\neq 0$, $T (m)$ and $V (m)$
are complex
vector bundles on which $z\in \T$ acts by fiberwise by multiplication
by $z^{m}$. Since the real
representations $\ur{\crep{m}}$ and $\ur{\crep{-m}}$ of $\T$ are
isomorphic, we are free to choose the signs of the $m_{j}$ and 
$m_{j}'$ which are not congruent to $0$ or $\nh$ modulo $n$ so that there
are integers $\ell_{j},r_{j}, \ell_{j}',$ and $r_{j}'$ satisfying
\begin{equation} \label{eq-ellj-rj}
\begin{split}
    m_{j}  & = n \ell_{j} + r_{j} \\
    m_{j}' & = n \ell_{j}' + r_{j}'
\end{split}
\end{equation}
with $0<r_{j}, r_{j}' < \nh$.    

We must choose orientations for $T (0)$, $T_{0}$, and $T_{\nh}$;  our 
choice will depend on the parity of $n$.

If $n$ is odd, then the complex structures on $T_{0<r<\nh}$ and the
orientation of $T$ itself determine an orientation for $T_{0}$.  We
choose the orientation of $T (0)$ and the signs of the $m_{j}\equiv 0$
so that \eqref{eq-T-zero-T-sub-zero} respects the orientations on both
sides.  In this case let $\signs = 0.$

If $n$ is even, then we simply choose an orientation for $T_{0}$; this
induces one on $T_{\nh}$ so that \eqref{eq-T-T-sub-zero-to-h} respects
the orientations.  We choose the orientation of $T (0)$ and the 
signs of the $m_{j}\equiv 0, \nh$ so that 
\[
    (T_{0}+T_{\nh})\restr{F} \cong  T (0) \oplus
\bigoplus_{m_{j}\equiv 0, \nh} \ur{T (m_{j})}.
\]
is compatible with the orientations on both sides.  For those $j$ such
that $m_{j}\equiv \nh$, we define $\ell_{j}$ so that
\[
      m_{j} = n \ell_{j} + \nh.
\]
We have two possibly distinct orientations on $T_{0}\restr{F}$:
one is obtained by restricting our chosen orientation; the other is
induced by the the decomposition \eqref{eq-T-zero-T-sub-zero}.
We define 
\[
       \signs\: \pi_{0} (X^{\T}\cap P) \xra{} \{0,1\}
\]
to be zero if these orientations agree, and one if they do not.  Note
that the decomposition \eqref{eq-T-h-decomp} of $T_{\nh}\restr{F}$
respects the 
orientations on both sides if and only if $\signs (F) =0$.

We choose the orientations and rotation numbers for $V$ in the same
way, and so define a function $\signs'$.  

Let us write 
\begin{align*}
      e_{r} & = \rank_{\C} T_{r} && 0<r<\nh \\
      e_{r}'& = \rank_{\C} V_{r} && 0<r<\nh \\
      e_{r} & = \thalf \rank_{\R} T_{r} && r = 0,\nh \\
      e_{r}'& = \thalf \rank_{\R} V_{r} && r = 0, \nh \\
      d_{j} & = \rank_{\C} T (m_{j}) \\
      d_{j}'& = \rank_{\C} V (m_{j}').
\end{align*}
Note that $e_{r}$ is an integer if $P$ contains an $F$ which
is non-empty.
For $m_{j}\neq 0$ let $x_{j,1},\dots ,x_{j,d_{j}}$ be the roots of the
total Chern class 
of $T (m_{j})$, and define $x_{j,1}',\dots ,x_{j,d_{j}'}'$ to be the
roots of the total Chern class of $V (m_{j}')$.

In terms of these we define several quantities which appear repeatedly
in our analysis. 
\begin{align*}
      \epsilon & = \lch{\lambda (\signs' + \sum d_{j}' \ell_{j}'
                   -\signs - \sum d_{j}\ell_{j})} \\
      \alpha   & = -\frac{1}{2n} \sum_{r=1}^{\nh} (e_{r}' - e_{r}) r^{2} \\
      G        & = \sum_{j} \frac{n}{2}  ( d_{j}'{\ell_{j}'}^{2} -
d_{j}\ell_{j}^{2})  +  d_{j}' \ell_{j}' r_{j}' - d_{j}\ell_{j} r_{j} \\
      H        & = \sum_{j} \ell_{j}' \borel{c_{1} V (m_{j}')} 
                -  \sum_{j} \ell_{j}  \borel{c_{1}T (m_{j})} \\
               & = \sum_{j} \ell_{j}' ( d_{j}' m_{j}'z  + \sum_{i} x_{j,i}') 
           -  \sum_{j} \ell_{j}  ( d_{j}  m_{j} z  + \sum_{i} x_{j,i}).
\end{align*}

Let $\Nu$ be the complex $\T$-line bundle 
\[
   \Nu = \otimes_{0<r<\nh} \det (V_{r})^{-r} \det (T_{r})^{r}
\]
over $X^{\T[n]}$.

We recall that $a$ is a special point of exact order $n$.  Let
$\lambda$ be the lattice point $na$.  Let $\gamma=\gamma (\lambda)$ be 
the constant in the translation formula
\eqref{eq-theta-translates-by-char} for $\theta$, and let $S (x) =
e^{\gamma x}$.  

Our hypotheses imply the following results.      They are essentially
contained in  \cite{BottTaubes:Rig}, to which we refer for the proof
of the first lemma.  The next three lemmas are trivial in the case
$\gamma=0$; proofs in the case that \eqref{eq-ccr} holds will be given
in \S\ref{sec-ccr-consequences}, below.  

\begin{Lemma} \label{t-le-signs-odd}
The quantity $\epsilon$
is independent of the component $F$ of $X^{\T}$, provided that $F$ is
contained in $P$.  In other words, $\epsilon$ is a locally constant
function on $X^{\T[n]}.$
\end{Lemma}

\begin{proof}
This is Lemma 9.3 in \cite{BottTaubes:Rig}.   
\end{proof}

\begin{Lemma} \label{t-le-alpha}
We have an equality 
\[
 S (\alpha) =  S ( G);
\]
in particular, this quantity is constant on $P$.
\end{Lemma}

\begin{Lemma} \label{t-le-beta-odd}
If $n$ is odd, then the complex $\T$-line bundle $\Nu$ has an $n$\th
root $\Nu^{1/n}$ over $P$, with the property that
\[
       S (i^{*}\borel{c_{1} ( \Nu^{1/n})}) =  S (H).
\]
\end{Lemma}

\begin{Lemma} \label{t-le-beta-even}
If $n$ is even, then the complex $\T$-line bundle $\Nu$ has an $\nh$\th
root $\Nu^{\nhr}$ over $P$, with the property that
\[
      S (\tfrac{1}{2} i^{*}\borel{c_{1} (\Nu^{\nhr})} )
       =   S ( H + \tfrac{1}{2}\borel{c_{1} (i^{*} (V_{h}- T_{h}))}).
\]
\end{Lemma}

\subsection{Construction of $\Theta$: ordinary points}
\label{sec:constr-thet-ordin}

First we give a formula for $\eqvtht_{b}$ when $b$ is ordinary.
Taking $a=0$  and $\eqvtht_{0}=1$ in \eqref{eq-mu-transition-fn}
gives the formula  
\begin{equation}\label{eq-mu-ordinary}
       \eqvtht_{b}   =     \trans_{b}^{*}  ( e (0,b) i^{*}\eqvtht_{0})
                     =     \trans_{b}^{*}   e (0,b).
\end{equation}
A priori, this formula determines $\eqvtht_{b}$ only on open subsets
$U$ of $U_{b}\backslash \{b \}$ and only on those $U_{b}$ such that $U_{b}\cap
U_{0}$ is non-empty, but in fact it determines $\eqvtht_{b}$ on all of
$U_{b}$ for any $b$.   
With our notations, we have 
\begin{equation}\label{eq-mu-ordinary-explicit}
   e (0,b) = \frac{\prod_{j} \prod_{i=1}^{d_{j}'}
                      \theta (x'_{j,i} + m_{j}'z )}
         {\prod_{j} \prod_{i=1}^{d_{j}}
                      \theta (x_{j,i} + m_{j}z )}.
\end{equation}
If $U_{b}\cap U_{0}$ is non-empty, or if $U_{b}\cap U_{a}$ is empty
for all special points $a$, then we orient $T (0)$ and $V (0)$ by
choosing the $m_{j}$ and $m_{j}'$ to be positive.  Otherwise, there is
a unique $a$ such that $U_{a}\cap U_{b}$ is non-empty, and we follow
the procedure described above to orient $T (0)$ and choose signs for
the $m_{j}$.
 
\begin{Lemma} \label{t-le-e-b-elliptic}
As a function of $z$ we have 
\[
   e (0,b) (z+\lm) = e (0,b) (z).
\]
In other words, the formula \eqref{eq-mu-ordinary-explicit} defines a
global section of $\HT (X^{b})\otimes \mero{C}$.  
\end{Lemma}

The proof is given in \S\ref{sec-ccr-consequences}.

\begin{Lemma} \label{t-le-theta-b-done}
The formula \eqref{eq-mu-ordinary} defines
an element of $\ET (X)_{b} (U_{b})$.
\end{Lemma}

\begin{proof}
What must be shown is that $\trans_{b}^{*}e (0,b)$ has no pole at
the origin.  Suppose it does; we shall show that $b$ is a special point.
Choose a lift of $b$ to $\C$; we may call it $b$ in view of Lemma
\ref{t-le-e-b-elliptic}.  We have  
\[
   \trans_{b}^{*}e (0,b) =     \frac{\prod_{j} \prod_{i=1}^{d_{j}'}
                      \theta (x'_{j,i} + m_{j}'z + m_{j}'b)}
         {\prod_{j} \prod_{i=1}^{d_{j}}
                      \theta (x_{j,i} + m_{j}z +m_{j}b)}.
\]
If this has a pole at the origin then $\theta$ has a zero at $m_{j}b$
for some $m_{j}$ or a pole at $m_{j}'b$ for some $m_{j}'$.  Let us
take the first case for definiteness.  By assumption, the zeros and
poles of $\theta$ are contained in $N^{-1}\Lambda$, so 
\[
    m_{j} b \in N^{-1}\Lambda
\]
and 
\[
   N m_{j} b\in \Lambda.
\]
By Definition~\ref{def-special-vector-bundle}, $b$ is a
special point.
\end{proof}

\subsection{Construction of $\Theta$: special points}
\label{sec:constr-thet-spec}

If $a$ is special, then as usual $U\subset U_{a}\cap U_{b}$ is
nonempty only if $b$ is ordinary.  Combining the two equations
\eqref{eq-mu-transition-fn} and \eqref{eq-mu-ordinary} gives
\begin{equation}
    \trans_{b}^{*}  e (0,b) = \trans_{b-a}^{*} ( e(a,b) i^{*}\eqvtht_{a})
\end{equation}
or equivalently 
\begin{equation}\label{eq-transfer}
     e (a,b)^{-1} \trans_{a}^{*}e (0,b) = i^{*}\eqvtht_{a}.
\end{equation}

\subsubsection{Special points of odd order}

In this section we give the formula for $\eqvtht_{a}$, supposing that $n$
is odd.

For $0<r<h$, let $Q_{r}$ be the exponential
characteristic class for complex vector bundles defined by the formula 
\[
    Q_{r} (L) = \theta (c_{1}L + r a).
\]

Let  $\eqvtht_{a} \in \ET (X) ( U_{a})$ be given by the formula 
\begin{equation}\label{eq-mu-a-odd}
    \eqvtht_{a} = 
\epsilon S (a \alpha ) S (\borel{c_{1} (\Nu^{1/n})}) 
\prod_{0<r<\nh }\frac{Q_{r} (V_{r})}{Q_{r} (T_{r})}.
\end{equation}

\begin{Proposition} \label{t-pr-transfer-special-points-odd}
If $n$ is odd, then the class $\eqvtht_{a}$ in \eqref{eq-mu-a-odd} satisfies
the equation  \eqref{eq-transfer}.
\end{Proposition}

\begin{proof}
We have 
\[
   e (0,b) = 
    \frac{\prod_{j} \prod_{i=1}^{d_{j}'}
                      \theta (x'_{j,i} + m_{j}'z )}
         {\prod_{j} \prod_{i=1}^{d_{j}}
                      \theta (x_{j,i} + m_{j}z )}.
\]
The iterated transformation formula \eqref{eq-trans-form-sg-it} gives
\begin{align}
    \trans_{a}^{*}e (0,b)  = &  
    \frac{\prod_{j} \prod_{i=1}^{d_{j}'}
                      \theta (x'_{j,i} + m_{j}'z + m_{j}'a )}
         {\prod_{j} \prod_{i=1}^{d_{j}}
                      \theta (x_{j,i} + m_{j}z + m_{j}a )}  \notag \\ 
                           = & 
 \epsilon S (a G + H) \prod_{0\leq r< \nh} \Pi_{r},
    \label{eq-trans-e-odd}  
\end{align}
where
\[
    \Pi_{r} = 
      \frac{\prod_{m_j' \equiv r \pod{n}}\prod_{i=1}^{d_{j}'} 
                \theta (x'_{j,i} + m_{j}'z + ra)}
           {\prod_{m_j \equiv r \pod{n}}\prod_{i=1}^{d_{j}} 
                \theta (x_{j,i} + m_{j} z + r a)}
\]
for $0\leq r <\nh$.

In view of Lemmas \ref{t-le-alpha} and \ref{t-le-beta-odd}, it remains to
observe that 
\begin{equation}\label{eq-Pi-zero-odd}
\Pi_{0} = e (a,b),
\end{equation}
while for $0<r<\nh$ we have
\begin{equation}\label{eq-Pi-r-odd}
    \Pi_{r} = i^{*} 
\left(\frac{\borel{Q_{r} (V_{r})}}{\borel{Q_{r} (T_{r})}}\right).
\end{equation}
\end{proof}

\subsubsection{Special points of even order}

In this section we give the formula for $\eqvtht_{a}$, supposing that $n$
is even.

Once again, for $0<r<  \nh$, let $Q_{r}$ be the power series 
\[
    Q_{r} (x) = \theta (x + r a).
\]
For all $r$ these give characteristic classes for complex vector
bundles.  

Now, however, for $r=\nh$, set 
\[
      Q_{\nh} (x) = S (-\tfrac{1}{2}x) \theta (x + \lmh).
\]
The formulae
\begin{align*}
     \theta (-x) & = - \theta (x) \\
     \theta (x+\lambda) & = \lch{\lambda} S (x+\lmh) \theta (x)
\end{align*}
imply the following.

\begin{Lemma} \label{t-le-Qnh-odd}
The power series $Q_{\nh}$ satisfies 
\[
    Q_{\nh} (-x) = - \lch{\lm} Q_{\nh} (x),
\]
and so defines a characteristic class of oriented even real vector bundles. 
\qed
\end{Lemma}


Let  $\eqvtht_{a} \in \ET (X)_{a} (U_{a})$ be given by the formula 
\begin{equation}\label{eq-mu-a-even}
       \eqvtht_{a} = 
              \epsilon 
              S (a \alpha ) S (\thalf \borel{c_{1} (\Nu^{\nhr})}) 
              \prod_{0<r\leq \nh }\frac{Q_{r} (V_{r})}
                                       {Q_{r} (T_{r})}.
\end{equation}
We have used Lemma \ref{t-le-Qnh-odd} to ensure that the class
$Q_{\nh} (V_{\nh})/ Q_{\nh} (T_{\nh})$ is well-defined.
The proof of Theorem \ref{t-th-thom-theta} is completed by the
following.

\begin{Proposition} \label{t-pr-transfer-special-points-even}
If $n$ is even, then the class $\eqvtht_{a}$ in \eqref{eq-mu-a-even} satisfies
the equation  \eqref{eq-transfer}.
\end{Proposition}

\begin{proof}
The argument is similar to the proof in the odd case, Proposition
\ref{t-pr-transfer-special-points-odd}.  Once again we have 
\[
    \trans_{a}^{*}e (0,b)  =
     \lch{\lambda (\sum d_{j}' \ell_{j}' -\sum d_{j}\ell_{j})} 
     S (a G + H) \prod_{0\leq r\leq \nh} \Pi_{r}.
    \label{eq-trans-e-even}  
\]
The result follows from Lemmas \ref{t-le-alpha},
\ref{t-le-beta-even}, and the equations
\begin{align*}
\Pi_{0} & = (-1)^{\signs' - \signs} e (a,b) \\
\Pi_{r} & = i^{*} \left(\frac{\borel{Q_{r} (V_{r})}}{\borel{Q_{r}
(T_{r})}}\right) && 0<r<  \nh \\
\Pi_{h} & = (-c (\lm))^{\signs' - \signs}
S (\tfrac{1}{2} \borel{c_{1} (i^{*} ( V_{h} - T_{h}))})
i^{*} \left(\frac{\borel{Q_{h} (V_{h})}}{\borel{Q_{h}
(T_{h})}}\right).
\end{align*}
\end{proof}

\section{Equivariant elliptic cohomology of Thom spaces of virtual
representations} \label{sec:virtual-reps}

Let $V$ be a complex representation of $\T$,  considered as a
$\T$-vector bundle over a point.    The fixed sub-bundles $V^{\T[n]}$
are automatically orientable, so Lemma \ref{t-le-cocycle-defined} and
Proposition \ref{t-pr-thom-sheaf} apply even without the spin
hypothesis.  In particular, the elliptic cohomology $\ET (V)$ of the
Thom space of $V$, and more generally of a virtual complex
representation of $\T$, is an invertible sheaf of $\ET (\point) =
\O{C}$-modules, i.e. a holomorphic line bundle over $C$.  It is
illuminating to understand the statement of Theorem
\ref{t-th-thom-theta} in this case.

Recall that we have chosen a generator $z$ of the character group of
$\T$, and so identified the representation ring of $\T$ with a ring of
Laurent polynomials
\[
    R[\T] \cong \Z[z,z^{-1}].
\] 
If $f$ 
is a Laurent polynomial in $z$, we shall write $V (f)$ for the
associated virtual complex representation of $\T$. 

If $n$ is a natural number, let $C[n]$ be the subgroup of $C$
consisting of points of order $n$.  Let $\I (C[n]) \subset \O{C}$ be
the sheaf of ideals consisting of germs of holomorphic functions which
vanish at $C[n]$; it is a holomorphic line bundle over $C$.  If $C[n]$
is regarded as a divisor on $C$, then $\I (C[n])$ coincides with the
line bundle which is usually denoted $\O{C} (-C[n])$.

\begin{Proposition}\label{t-le-thom-sheaf-z-n}
The line bundles
$\ET (V)$ for $V \in R[\T]$ are determined up to isomorphism by  the
following. 
\begin{align}
\notag
 \ET (V (0)) & = \O{C} \\
\label{eq-Thom-sheaf-z-n}
 \ET (V (z^{n}))  & \cong \I (C[n])  \\
\intertext{and}
\notag
  \ET (V+W) & \cong \ET (V) \tensor{\O{C}} \ET (W).
\end{align}
for $V,W\in R[\T]$.
\end{Proposition}

\begin{proof}
The first part is clear, since $\ET (V (0))= \ET (\point)$.
The third part follows from the multiplicative property of the Thom
isomorphism \eqref{eq-thom-iso-mult} (and so of its associated euler
class), together with the construction \eqref{eq-twisted-sheaf} of
$\ET (X)^{[\phi,V]}\cong \ET (V)$.  

For \eqref{eq-Thom-sheaf-z-n}, let $V=V (z^{n})$, and consider the
multiplicative orientation $\phi$ given by the $\sigma$ function.  
A point $a\in C$ is special if and only if $na = 0$, and for such $a$
we have $V^{a} = V$.  For ordinary $b$ we have $V^{b} = 0$.  It
follows that 
\[
    e (a,b) = \sg (nz),
\]
a function which vanishes at precisely the points $w\in \C$ such that
$nw\in \Lambda$.  Thus in the gluing \eqref{eq-twisted-sheaf}, a
trivialization of $\ET (V)_{a}$ corresponds to a section of $\ET
(V)_{b}$ which vanishes at $a$; this is a description 
in terms of cocycles of the ideal sheaf $\I (C[n])$.  The isomorphism 
\[
   \ET (X)^{[\phi,V]}\cong \ET (V)
\]
of Proposition \ref{t-pr-thom-sheaf} gives the result.
\end{proof}

If $f = \sum d_{j} z^{m_{j}}$, let $D (f)$ be the divisor
\[
       D (f) = -\sum d_{j} C[m_{j}]
\]
(The subgroup $C[m_{j}]$ is considered as a divisor and the sum is as
divisors). 
Proposition \ref{t-le-thom-sheaf-z-n} says that there is an
isomorphism 
\[
    \ET (V (f))\cong \O{C} (D (f)).
\]
The degree of this line bundle is 
\begin{equation}\label{eq-deg-ET-f}
\deg \ET (V (f))  =  -\sum d_{j} m_{j}^{2}.
\end{equation}

\begin{Corollary} \label{t-co-degree-thom-sheaf}
If $V\in R[\T]$ is a complex virtual representation of $\T$, then 
\[
     \borel{p_{1} (V)} = -z^{2}\cdot \deg \ET (V) \in H^{4} (B\T;\Z).
\]\qed
\end{Corollary}

The application of Theorem \ref{t-th-thom-sigma} to complex
representations of $\T$ is as follows.  

\begin{Proposition} \label{t-pr-thom-sheaf-rep-trivial}
The line bundle $\ET (V)$ associated to a virtual representation $V$
of $\T$ is trivial precisely when $\borel{p_{1} (V)} = 0.$  If 
\[
     f  = \sum d_{j} z^{m_{j}}
\]
and $\borel{p_{1} (V (f))} = 0$, then the meromorphic function 
\begin{equation}\label{eq-sigma-of-virt-rep}
      \prod_{j} \sg (m_{j} z)^{d_{j}}
\end{equation}
defines a holomorphic trivialization of $\O{C} (D(f))\cong \ET (V
(f))$.  If in addition $\borel{w_{2} (V)} = 0,$ it is precisely the
trivialization given by Theorem \ref{t-th-thom-sigma}. 
\end{Proposition}

\begin{proof}
Let 
\[
D  = \sum n_{P} (P)
\]
be a divisor on $C$.  Recall that the line bundle $\O{C} (D)$
associated to $D$ is trivial precisely when
\begin{subequations}
\begin{gather}\label{eq-degree-condition}
    \deg D = \sum n_{P} = 0 \\
    \sideset{}{^{C}}\sum n_{P} (P) = 0; \label{eq-group-condition}
\end{gather}
\end{subequations}
the second sum is taken in the group structure of the elliptic curve $C$.
For the divisors we are considering, the condition
\eqref{eq-group-condition} is always satisfied, since as a group
\[
    C[n] \cong (Z/n)^{2},
\]
and 
\[
    \sum_{g\in (Z/n)^{2}} g = 0.
\]
It follows that $\ET (V)$ is trivial precisely when $\deg \ET (V)
= 0$, that is when $\borel{p_{1} (V)} =0$.

It is well-known and easy to check that the product of sigma functions
\eqref{eq-sigma-of-virt-rep} descends to a trivialization of $\O{C} (D
(f))$, and it is also easy to check that it coincides with the
trivialization of Theorem \ref{t-th-thom-sigma}.   
\end{proof}

Proposition \ref{t-pr-thom-sheaf-rep-trivial} raises the question of
the role of the condition \eqref{eq-swcr} on the equivariant second
Stiefel-Whitney class.  One way to say this is as follows.  Let 
\[
f =  \sum d_{j} z^{m_{j}}
\]
as in the Proposition, and let $g$ be the associated trivialization of
$\O{C} (D (f))$: 
as a function on $\C$,
\[
g (z)    = \prod_{j} \sg (m_{j} z)^{d_{j}}.
\]
Let $\iota: C\to C$ be the involution 
\[
    \iota (P) = -P.
\]
The equality of divisors 
\[
    \iota^{*} D (f) = D (f)
\]
gives a canonical isomorphism 
\[
    \iota^{*}\O{C} (D (f)) \cong \O{C} (D (f)), 
\]
of line bundles over $C$, and it is natural to ask whether 
\begin{equation}\label{eq-g-real}
      \iota^{*} g = g
\end{equation}
under this isomorphism.  As a function on $\C$, this amounts to asking whether
\[
     g (-z) = g (z).
\]

The function $g$ is uniquely determined among trivializations of
$\O{C} (D (f))$ by the equation 
\[
   g' (0) = \sum d_{j} m_{j},
\]
and the parity of this quantity is the second Stiefel-Whitney class
\[
    \borel{w_{2} (V)} = g' (0) z \in H^{2} (B\T;\Z/2).
\]
Thus
\[
      \iota^{*}g = g
\]
precisely when $\borel{w_{2} (V)} = 0$.  This means, for example, that
the trivialization $g$ is independent of the choice of generator $z$
of the character group of $\T$.

\section{Consequences of the characteristic class restrictions}

\label{sec-ccr-consequences}

In this section we prove Lemmas \ref{t-le-alpha}, \ref{t-le-beta-odd},
and \ref{t-le-beta-even}.
We retain the set-up and notations of the previous section.

First let us work out explicitly some consequences of the
characteristic class restrictions \eqref{eq-ccr}.

\begin{Lemma} \label{t-le-ccr-imply-always}
If equations \eqref{eq-ccr} hold, then we have 
\begin{align}
\label{eq-ccr-z-3}
 \sum_{j} d_{j} m_{j} & \equiv \sum_{j} d_{j}' m_{j}' \mod{2} \\
\label{eq-ccr-z-2}
 \sum_{j} d_{j} m_{j}^{2} & = \sum_{j} d_{j}' {m_{j}'}^{2} \\
\label{eq-ccr-z-1}
 \sum_{j} m_{j} \sum_{i=1}^{d_{j}} x_{j,i} & = 
 \sum_{j} m_{j}' \sum_{i=1}^{d_{j}' }x_{j,i}'.
\end{align}
\end{Lemma}

\begin{proof}
Introduce formal roots $y_{0,j}$ and $y_{0,j}'$ so that 
\begin{align*}
     1 + p_{1} T (0) + \dots  & = \prod (1 - y_{0,i}^{2}) \\
     1 + p_{1} V (0) + \dots  & = \prod (1 - {y_{0,i}'}^{2}).
\end{align*}
The equation \eqref{eq-phcr} 
\[
   \borel{\phalf (V - T)}  = 0
\]
implies that 
\[
   \borel{\phalf (V\restr{M^{\T}} - T\restr{M^{\T}})} = 0.
\]
On the other hand, this class is given by half the degree-four component of 
\[
 \frac{\prod (1 - {y_{0,i}'}^{2})
 \prod_{j}
\prod_{i=1}^{d_{j}'} 
 (1- ({x_{j,i}'}^{2} +2 m_{j}'x_{j,i}'z + {m_{j}'}^{2}z^{2}))}
      {\prod (1 - y_{0,i}^{2})
\prod_{j}
\prod_{i=1}^{d_{j}} 
 (1- (x_{j,i}^{2} +2 m_{j}x_{j,i}z + m_{j}^{2}z^{2}))}
\]
Examining the coefficient of $z$ gives \eqref{eq-ccr-z-1}, and
examining the coefficient of $z^{2}$ gives \eqref{eq-ccr-z-2}.
Equation \eqref{eq-ccr-z-3} follows from the equation $\borel{w_{2}
(T-V)} = 0$ by a similar argument.
\end{proof}

\begin{Lemma} \label{t-le-ccr-imply-odd-even}
If equations \eqref{eq-ccr} hold, and $n$ is odd, then
\begin{equation}
\label{eq-ccr-z-n-odd}
 \sum_{0<r<\nh} r (c_{1} T_{r} - c_{1} V_{r})  \equiv 0 \pod{n}.
\end{equation}
If $n$ is even, then 
\begin{equation}\label{eq-ccr-z-n-even}
\begin{split} 
 \sum_{0<r<\nh} r (c_{1} T_{r} - c_{1} V_{r}) & \equiv 0 \pod{\nh} \\
 \frac{1}{h}\left( \sum_{0<r<\nh} r (c_{1} T_{r} - c_{1} V_{r})
 \right)_{\mod{2}}  
   & = 
 w_{2} (V_{\nh} - T_{\nh}).
\end{split}
\end{equation}
\end{Lemma}

\begin{proof}
We treat the case that $n$ is even; the case that $n$ is odd is
similar.  Let $z_{n}=z\restr{B\T[n]}$.  Introduce formal
roots $y_{r,j}$ and $y_{r,j}'$ so that 
\begin{align*}
   1 + c_{1} T_{r} + \dots  & = \prod ( 1 +  y_{r,j})    && 0<r<\nh \\
   1 + c_{1} V_{r} + \dots  & = \prod ( 1 + y_{r,j}')    && 0<r<\nh \\
   1 + p_{1} T_{r} + \dots  & = \prod ( 1 - y_{r,j}^{2}) && r = 0, \nh \\
   1 + p_{1} V_{r} + \dots  & = \prod ( 1 - {y_{r,j}'}^{2}) && r = 0,\nh.
\end{align*}
Then $\phalf (\borel{T} - \borel{V})\restr{B\T[n]\times M^{\T[n]}}$
is given by minus half the degree-four component of 
\[
 \frac{\prod_{0\leq r\leq \nh} 
       \prod_{i=1}^{e_{r}} 
 (1- (y_{r,i}^{2} +2 r y_{r,i}z_{n} + r^{2}z_{n}^{2}))}
      {\prod_{0\leq r\leq \nh} 
       \prod_{i=1}^{e_{r}'} 
 (1- ({y_{r,i}'}^{2} +2 r y_{r,i}'z_{n} + r^{2}z_{n}^{2}))}.
\]
The coefficient of $z_{n}$ is 
\[
   \sum_{0<r<\nh} r (c_{1} T_{r} - c_{1} V_{r}) + \nh \sum_{j}
   y_{\nh,j}  - y_{\nh,j}'. 
\]
The characteristic class restriction \eqref{eq-ccr}  implies that this
quantity is zero; the claims of the lemma follow.
\end{proof}

\begin{proof}[Proof of Lemma \ref{t-le-e-b-elliptic}]
We have 
\begin{align*}
   e (0,b) (z + \lm)  & = \frac{\prod_{j} \prod_{i=1}^{d_{j}'}
                      \theta (x'_{j,i} + m_{j}'z + m_{j}'\lm)}
         {\prod_{j} \prod_{i=1}^{d_{j}}
                      \theta (x_{j,i} + m_{j}z +m_{j}\lm)} \\
                    &  = 
        \lch{\lm (\sum d_{j}' m_{j}' - d_{j}m_{j})} \\
&\qquad        S (\sum m_{j}' (d_{j}' m_{j}' z + \sum x_{j,i}') - 
           \sum m_{j} (d_{j} m_{j} z + \sum x_{j,i})) \\
&\qquad       S (\lmh (\sum d_{j}'  {m_{j}'}^{2} - \sum d_{j} m_{j}^{2}))\\
 &\qquad           e (0,b) (z) \\
 &                     =  e (0,b) (z).
\end{align*}
The third equation uses Lemma \ref{t-le-ccr-imply-always}.
\end{proof}

\begin{proof}[Proof of Lemma \ref{t-le-alpha}]
Recall that
\[
G  = \sum_{j} \frac{n}{2}  ( d_{j}'{\ell_{j}'}^{2} -
d_{j}\ell_{j}^{2})  +  d_{j}' \ell_{j}' r_{j}' - d_{j}\ell_{j} r_{j}.
\]
Substituting the equation 
\[
     m_{j} = n \ell_{j} + r_{j}
\]
(and similarly for $m_{j}'$) into equation \eqref{eq-ccr-z-2} gives 
\begin{align*}
   0 & = \sum_{j} d_{j}' {m_{j}'}^{2} - d_{j} m_{j}^{2} \\
     & = 
\sum_{j} 
 d_{j}' (n^{2} {\ell_{j}'}^{2} + 2 n \ell_{j}'r_{j}' +  {r_{j}'}^{2}) -
\sum_{j} 
 d_{j} (n^{2} {\ell_{j}}^{2} + 2 n \ell_{j}r_{j} + r_{j}^{2}).
\end{align*}
It follows that 
\[
2 n G = -  \sum_{r} (e_{r}' - e_{r}) r^{2},
\]
and that this quantity is divisible by $n$.
\end{proof}

\begin{proof}[Proof of Lemma \ref{t-le-beta-odd}]
The first Chern class of the bundle
\[
   \Nu = \prod_{0<r<\nh} \det (V_{r})^{-r} \det (T_{r})^{r}
\]
is
\[
   c_{1}\Nu = -\sum_{r} r (c_{1}V_{r} - c_{1} T_{r}).
\]
The restriction of $\Nu$ to $M^{\T}$ has Chern class
\[
    \borel{c_{1} (\Nu\restr{M^{\T}})}  = 
    -\sum_{j} r_{j}' (d_{j}' m_{j}' z + \sum_{i=1}^{d_{j}'} x_{j,i}') +
    \sum_{j} r_{j}  (d_{j} m_{j} z + \sum_{i=1}^{d_{j}} x_{j,i}).
\]
Adding to this zero in the form (from \eqref{eq-ccr-z-1} and
\eqref{eq-ccr-z-2})  
\[
    \sum_{j} m_{j}' (d_{j}' m_{j}' +  \sum_{i=1}^{d_{j}' }x_{j,i}')
 -  \sum_{j} m_{j} (d_{j} m_{j} + \sum_{i=1}^{d_{j}} x_{j,i})
\]
gives 
\begin{align*}
\borel{c_{1} (\Nu\restr{M^{\T}})}  & = 
\sum_{j} (m_{j}'- r_{j}') (d_{j}' m_{j}' z + \sum_{i=1}^{d_{j}'} x_{j,i}')- 
\sum_{j} (m_{j} - r_{j} ) (d_{j} m_{j} z + \sum_{i=1}^{d_{j}} x_{j,i}) \\
    & = 
n \sum_{j} \ell_{j}' (d_{j}' m_{j}' z +  \sum_{i=1}^{d_{j}'} x_{j,i}') -
n \sum_{j} \ell_{j}  (d_{j}' m_{j} z + \sum_{i=1}^{d_{j}}   x_{j,i})
\\
   & = n H.
\end{align*}
\end{proof}

\begin{proof}[Proof of Lemma \ref{t-le-beta-even}]
First observe that if 
\[
   \Nu = \prod_{0<r<\nh} \det (V_{r})^{-r} \det (T_{r})^{r},
\]
then Lemma \ref{t-le-ccr-imply-odd-even} shows that $\Nu$ has an $\nh$\th
root, with the property that
\[
    c_{1}\Nu^{\nhr}_{\mod{2}} = w_{2} (V_{\nh} - T_{\nh}).
\]

The restriction of $\Nu$ to $M^{\T}$ has Chern class
\[
    \borel{c_{1} (\Nu\restr{M^{\T}})}  = 
-\sum_{0<r_{j}'<\nh} r_{j}' (d_{j}' m_{j}' z + \sum_{i=1}^{d_{j}'} x_{j,i}') +
 \sum_{0<r_{j}<\nh} r_{j}  (d_{j} m_{j} z + \sum_{i=1}^{d_{j}} x_{j,i}).
\]
Adding to this zero in the form (from \eqref{eq-ccr-z-1} and
\eqref{eq-ccr-z-2})  
\[
    \sum_{j} m_{j}' (d_{j}' m_{j}'z +  \sum_{i=1}^{d_{j}' }x_{j,i}')
 -  \sum_{j} m_{j} (d_{j} m_{j} z+ \sum_{i=1}^{d_{j}} x_{j,i})
\]
gives 
\begin{align*}
\borel{c_{1} (\Nu\restr{M^{\T}})}  = & n H \\
&+\sum_{m_{j}'\equiv\nh} m_{j}' (d_{j}' m_{j}'z +  \sum_{i=1}^{d_{j}'}x_{j,i}')
-\sum_{m_{j}'\equiv\nh}n \ell_{j}' ( d_{j}' m_{j}'z  + \sum_{i} x_{j,i}') \\
&-\sum_{m_{j}\equiv\nh} m_{j} (d_{j} m_{j}z + \sum_{i=1}^{d_{j}} x_{j,i}) 
+\sum_{m_{j}\equiv\nh} n \ell_{j}  ( d_{j}  m_{j} z  + \sum_{i} x_{j,i}) \\
= & n H \\
& +
\nh \sum_{m_{j}'\equiv\nh} (d_{j}' m_{j}'z +  \sum_{i=1}^{d_{j}'}x_{j,i}') - 
\nh \sum_{m_{j} \equiv\nh} (d_{j} m_{j}z +  \sum_{i=1}^{d_{j}}x_{j,i}).
\end{align*}
\end{proof}







\end{document}